\documentclass{amsart}

\usepackage{amsmath, amsfonts, mathrsfs}
\usepackage{amssymb,latexsym, fontenc}
\usepackage{enumerate}
\usepackage{hyperref}
\usepackage{bbm}

\usepackage[top=3.1cm, bottom=3.1cm, left=3.1cm, right=3.1cm, includefoot, heightrounded]{geometry}

\newtheorem{theorem}{Theorem}[section]
\newtheorem{corollary}[theorem]{Corollary}
\newtheorem{lemma}[theorem]{Lemma}

\newtheorem{example}[theorem]{Example}

\theoremstyle{definition}
\newtheorem{definition}[theorem]{Definition}

\numberwithin{equation}{section}

\begin{document}

\title[Characterizations of Pettis Integrable Multifunctions]
{Characterizations of Pettis Integrable Multifunctions in locally convex spaces}

\author{ Sokol Kaliaj }

\address{
Mathematics Department, 
Science Natural Faculty, 
University of Elbasan,
Elbasan, 
Albania.
}

\email{printimesokol@gmail.com}

\thanks{}



\begin{abstract}
In this paper, we present characterizations of Pettis integrability in $bcc(X)$ and $cwk(X)$ 
of a multifunction $\Gamma : \Omega \to cc(X)$,  
where $(\Omega, \Sigma, \mu)$ is a complete finite measure space 
and 
$cc(X)$ $(bcc(X), cwk(X))$ is the family of all nonempty 
convex closed (bounded convex closed, convex weakly compact) subsets of a 
Hausdorff locally convex topological vector space $X$. 
Incidentally, a convergence theorem for scalarly integrable multifunctions is proved. 
As an application, 
we are going to prove full characterizations 
of Pettis integrability for fuzzy maps defined on $\Omega$ and taking values in the generalized fuzzy number hyperspace 
$\mathcal{F}_{cwk}(X)$. 
\end{abstract}

\maketitle

\section{Introduction}

There are well-known characterizations of Pettis integrability of a multifunction in non-separable Banch spaces. 
B.~Cascales, V.~Kadets and J.~Rodriguez in \cite[Lemma 2.1]{CASC_SelectorsA} 
have proved that 
a scalarly integrable multifunction $\Gamma : \Omega \to cwk(X)$ 
is Pettis integrable in $cwk(X)$ if and only if 
the map 
$\varphi_{\Gamma}^{A}: X^{*} \to \mathbb{R}$ defined by 
$\varphi_{\Gamma}^{A}(x^{*}) = \int_{A} \sigma(x^{*}, \Gamma(t)) d \mu$ for all $x^{*} \in X^{*}$ and $A \in \Sigma$  
is $\tau(X^{*},X)$-continuous, 
where $X^{*}$ is the topological dual of a Banach space $X$, 
$\sigma(\cdot, \Gamma(t))$ is the support function of  $\Gamma(t)$ 
and $\tau(X^{*},X)$ is the Mackey topology on $X^{*}$. 
K.~Musial has also proved in \cite[Theorem 1.4]{MUS_Pettis_MultiF} 
that if a multifunction $\Gamma : \Omega \to cc(X)$ is scalarly integrable 
and a map $\varphi_{\Gamma}: X^{*} \to L_{1}(\mu)$ is defined by 
$\varphi_{\Gamma}(x^{*}) = \sigma(x^{*}, \Gamma(\cdot))$ for all $x^{*} \in X^{*}$, then
\begin{itemize}
\item
$\Gamma$ is Pettis integrable in $cwk(X)$ 
if and only if 
the map $\varphi_{\Gamma}$ is $\tau(X^{*},X)$-$\sigma(L_{1},L_{1}^{*})$-continuous, 
where $\sigma(L_{1},L_{1}^{*})$ is the weak topology in $L_{1}(\mu)$,
\item
if $\Gamma$ is non-negative and Pettis integrable in $cwk(X)$, then $\varphi_{\Gamma}$ is $\tau(X^{*},X)$-$||\cdot||_{1}$-continuous, 
where $||\cdot||_{1}$ is the norm in $L_{1}(\mu)$.
\end{itemize}

In this paper, 
we show two new characterizations of Pettis integrability of a scalarly integrable multifunction $\Gamma : \Omega \to cc(X)$, 
see Lemma \ref{L_Pettis_BCC} and Theorem \ref{T_Pettis_Chr_CWK}. 
The first one shows that if 
$\Gamma$ is Pettis integrable in $bcc(X)$, 
then $\varphi_{\Gamma}$ is $\beta(X^{*}, X)$-$||\cdot||_{1}$-continuous, 
where $\beta(X^{*},X)$ is the strong topology on $X^{*}$.
The second one shows that $\Gamma$ is Pettis integrable in $cwk(X)$ if and only if 
$\varphi_{\Gamma}$  is  $\tau(X^{*}, X)$-$||\cdot||_{1}$-continuous.
To prove these results, we need to characterize the continuity of a sublinear map  
$\varphi : X^{*} \to L_{1}(\mu)$ ($\varphi : X^{*} \to \mathbb{R}$), 
see Lemma \ref{L_SublinearL1_Cont_Zero} and Lemma \ref{L_SublinearL1_Cont_Bounded}. 
Incidentally, a convergence theorem for scalarly integrable multifunctions is proved, see  
Theorem \ref{T_Convergence_Theorem}. 
Corollary \ref{C_Dominated_Convergence_Theorem} shows a version of the Lebesgue dominated convergence theorem 
for scalarly integrable multifunctions.

As an application of Theorem \ref{T_Pettis_Chr_CWK}, 
we are going to prove full characterizations of Pettis integrability 
for a fuzzy map $\Gamma : \Omega \to \mathcal{F}_{cwk}(X)$ 
in terms of the map
\begin{equation*}
\begin{split}
\varphi_{\Gamma} : X^{*} \to \prod_{\alpha \in (0,1]} L_{1}(\mu),
\quad
\varphi_{\Gamma}(x^{*}) = (\sigma(x^{*}, \Gamma^{\alpha}))_{\alpha \in (0,1]}, 
\end{split}
\end{equation*}
where $\Gamma^{\alpha}$ is the $\alpha$-cut of $\Gamma$, 
i.e., 
\begin{equation*}
\begin{split}
\Gamma^{\alpha}(t) = \{x \in X : \Gamma(t) \geq \alpha \}, 
\qquad
(t \in \Omega, \alpha \in (0,1]). 
\end{split}
\end{equation*}

The paper is organized as follows by sections. 
\begin{itemize} 
\item[(\S 2)]
We recall some basic concepts and facts such as Mackey topology $\tau(X^{*},X)$, 
strong topology $\beta(X^{*},X)$,  
multimeasures, Pettis integrable multifunctions, sublinear maps and a useful theorem for our purposes, 
Theorem \ref{T_CWK}.
\item[(\S 3)]
The continuity on $X^{*}$ of a 
sublinear map $\varphi: X^{*} \to L_{1}(\mu)$ ($\varphi: X^{*} \to \mathbb{R}$) 
is characterized by its continuity at $0 \in X^{*}$ 
under a  Hausdorff locally convex topology $\mathfrak{T}^{*}$ 
on the topological dual $X^{*}$ of $X$,  
see Lemma \ref{L_SublinearL1_Cont_Zero} (Lemma \ref{L_SublinearR_Cont_Zero}).
We introduce the notion of bounded sublinear maps (Definition \ref{D_Bounded_Sublinear}), 
and show that the class of continuous sublinear maps 
coincides with 
the class of bounded sublinear maps, 
see Lemma \ref{L_SublinearL1_Cont_Bounded} (Lemma \ref{L_SublinearR_Cont_Bounded}). 
\item[(\S 4)]
We first prove that if a multifunction 
$\Gamma : \Omega \to cc(X)$ is Pettis integrable in $bcc(X)$, 
then  $\varphi_{\Gamma}$ is $\beta(X^{*}, X)$-$||\cdot||_{1}$-continuous, 
see Lemma \ref{L_Pettis_BCC}. 
Then, Theorem \ref{T_Pettis_Chr_CWK} shows that 
a scalarly integrable multifunction $\Gamma : \Omega \to cc(X)$ is Pettis integrable in $cwk(X)$ 
if and only if  
$\varphi_{\Gamma}$ is  $\tau(X^{*}, X)$-$||\cdot||_{1}$-continuous.    
Theorem \ref{T_Pettis_Chr_CWK} is the key to prove Theorem \ref{T_Convergence_Theorem}. 
Corollary \ref{C_Dominated_Convergence_Theorem} shows a version of the Lebesgue dominated convergence theorem 
for scalarly integrable multifunctions. 
\item[(\S 5)] 
As an application of Theorem \ref{T_Pettis_Chr_CWK}, 
we will prove full characterizations of Pettis integrability 
for a fuzzy map, see Theorem \ref{T_FPettis_Chr_CWK}.
\end{itemize}

\section{Basic concepts and facts}

Throughout this paper, 
$X$ is a non-separable complete Hausdorff locally convex topological vector space 
with the topological dual $X^{*}$, and 
$(\Omega, \Sigma, \mu)$ is a complete finite measure space. 
We denote by $\mathcal{P}(X)$ the family of all nonempty subsets of $X$ 
and by $cc(X)$ ($bcc(X)$, $cwk(X)$) is denoted the subfamily of $\mathcal{P}(X)$ of 
all convex closed (bounded  convex closed,  convex weakly compact) subsets of $X$. 
For each $C \in \mathcal{P}(X)$ 
the \textit{support} of $C$ is defined by equality
\begin{equation*}
\begin{split}
\sigma(x^{*}, C) = 
\sup \{ \langle x^{*}, x\rangle : x \in C \}
\qquad (x^{*} \in X^{*}).
\end{split}
\end{equation*}
By $\langle X, X^{*}\rangle$ is denoted the ordinary dual system
(cf. \cite[p.123]{SCH}) with bilinear form defined by
\begin{equation*}
\begin{split}
\langle x, x^{*} \rangle =x^{*}(x)
\qquad
(x \in X, x^{*} \in X^{*}),
\end{split}
\end{equation*} 
$\tau(X, X^{*})$ is the \textit{Mackey topology} on $X$, 
i.e., the topology of uniform convergence on the family $\mathcal{M}(X^{*},X)$ of all absolutely convex (circled and convex) 
and $\sigma(X^{*}, X)$-compact subsets of $X^{*}$, cf. \cite[Corollary 1, p.131]{SCH}, 
$\tau(X^{*},X)$ is the \textit{Mackey topology} on $X^{*}$, 
i.e., 
the topology of uniform convergence on the family $\mathcal{M}(X,X^{*})$ of all absolutely convex 
and $\sigma(X, X^{*})$-compact subsets of $X$.
By $\beta(X^{*},X)$ is denoted the \textit{strong topology} in $X^{*}$, 
i.e., the topology of uniform convergence in $X^{*}$ 
on the family of all bounded subsets of $X$.
Thus, if $B$ is a nonempty bounded subset of $X$, 
then the absolute polar 
$B^{o} = \{x^{*} \in X^{*}: \sup_{x\in B} |x^{*}(x)| \leq 1\}$ of $B$ is a $0$-neighborhood 
under the $\beta(X^{*},X)$-topology and the function
\begin{equation*}
\begin{split}
p_{B} : X^{*} \to [0,+\infty),
\quad
p_{B}(x^{*}) = \sup_{x \in B} |x^{*}(x)|
\end{split}
\end{equation*} 
is the Minkowski functional of $B^{o}$, cf. \cite[p.39]{SCH}.
Given a continuous seminorm $p$ on $X$ we write 
$U_{p} = \{x \in X : p(x) \leq 1\}$ and denote by $U_{p}^{o}$ the absolute polar of $U_{p}$, 
i.e., $U_{p}^{o} = \{x^{*} \in X^{*} : \sup_{x \in U_{p}} |x^{*}(x)| \leq 1\}$.


If  $(B_{n})$ is a sequence with terms in $\mathcal{P}(X)$, 
the symbol $\sum_{n=1}^{+\infty} B_{n}$ 
denotes a formal series. 
The series $\sum_{n=1}^{+\infty} B_{n}$ is said to be 
\textit{unconditionally convergent} 
provided that for every choice $b_{n} \in B_{n}, n \in \mathbb{N}$, 
the series $\sum_{n=1}^{+\infty} b_{n}$ is unconditionally convergent in $X$. 
In this case we define
\begin{equation}\label{eqD_Multimeasures}
\sum_{n=1}^{+\infty} B_{n} = 
\left \{
\sum_{n=1}^{+\infty} b_{n}: b_{n} \in B_{n} \text{ for all }n \in \mathbb{N}
\right \}.
\end{equation}

\begin{definition}\label{D_Multimeasures}
A map $M: \Sigma \to bcc(X)$ is called a \textit{multimeasure},  
if the map $\sigma(x^{*}, M(\cdot))$ is a countable additive real valued measure for every 
$x^{*} \in X^{*}$. 
The mapping $M$ is said to be a \textit{strong multimeasure}, if
\begin{itemize}
\item[(i)]
$M(\emptyset) = \{\theta\}$, 
\item[(ii)]
for each sequence $(A_{n})$ of pairwise disjoint members of $\Sigma$, 
we have $M(\bigcup_{n} A_{n}) = \sum_{n} M(A_{n})$, 
where  $\sum_{n} M(A_{n})$ is defined by \eqref{eqD_Multimeasures}.  
\end{itemize}
It is clear that each strong multimeasure is a multimeasure. 
If $M$ is a multimeasure taking values in the hyperspace $cwk(X)$, then the converse is also true, 
see \cite[Theorem 2.10]{KAL_REND}. 
\end{definition}


\begin{definition}\label{D_Pettis_Multi}
A multifunction $\Gamma: \Omega \to cc(X)$ is said to be 
\textit{scalarly measurable} (\textit{scalarly integrable}) 
if the map $\sigma(x^{*}, \Gamma(\cdot))$ is measurable (Lebesgue integrable) 
whenever $x^{*} \in X^{*}$. 
The multifunction $\Gamma$ is said to be
\textit{Pettis integrable} in $bcc(X)$ ($cwk(X)$),  
if 
\begin{itemize}
\item
$\Gamma$ is scalarly integrable,
\item
for each $A \in \Sigma$ there exists $M(A) \in bcc(X)$ ($M(A) \in cwk(X)$) such that
\begin{equation*}
\begin{split}
\sigma
\left ( 
x^{*}, M(A) 
\right )
=
\int_{A} \sigma(x^{*}, \Gamma(t)) d \mu
\quad
\text{for every }x^{*} \in X^{*}.
\end{split}
\end{equation*} 
\end{itemize}
We call $M(A)$ the \textit{Pettis integral} of $\Gamma$ over $A$ and 
set $(P)\int_{A} \Gamma(t) d \mu = M(A)$. 
It is clear that the map 
$M$ defined by $M(A) = (P)\int_{A} \Gamma(t) d \mu$ for all $A \in \Sigma$, 
is a multimeasure. 
\end{definition} 

The Pettis integral for multifunctions was first considered by Castaing and Valadier 
\cite[Chapter V]{CAST}. 
For the general theory of Pettis integrable multifunctions  
we use the papers 
of H.~Ziat, K.~El Amri and C.~Hess \cite{EL, ZIAT}, 
B.~Cascales, V.~Kadets, and J.~Rodr\'{i}guez \cite{CASC_SelectorsA, CASC_SelectorsB} and 
K.~Musial \cite{MUS_Pettis_MultiF, MUS_Approximation}. 
A presentation of the general theory of Pettis integration of vector functions can be
found in \cite{TAL} and \cite{MUS_Pettis_Integral_Functions}.

\begin{definition}\label{D_Uniformly_Integrable}
A family $\mathcal{F} \subset L_{1}(\mu)$ is said to be 
\textit{uniformly integrable} if $\mathcal{F}$ is bounded in $(L_{1}(\mu), ||\cdot||_{1})$ 
and for each $\varepsilon >0$ there exists $\delta >0$ such that 
for each $A \in \Sigma$, 
\begin{equation*}
\begin{split}
\mu(A) < \delta 
\Rightarrow
\sup_{f \in \mathcal{F}} \int_{A} |f(t)| d \mu < \varepsilon.
\end{split}
\end{equation*}
According to \cite[Corollary IV.8.10]{DUN} and \cite[Corollary IV.8.11]{DUN}, 
if $\mathcal{F} \subset L_{1}(\mu)$ is bounded in $L_{1}(\mu)$, then 
\begin{equation*}
\begin{split}
\mathcal{F}\text{ is uniformly integrable} \Leftrightarrow 
\mathcal{F}\text{ is weakly sequentiality compact},
\end{split}
\end{equation*}
and by Eberlein-$\check{S}$mulian theorem (\cite[Theorem V.6.1]{DUN}), 
it follows that 
\begin{equation*}
\begin{split}
\mathcal{F}\text{ is uniformly integrable} \Leftrightarrow 
\mathcal{F}\text{ is relatively weakly compact}.
\end{split}
\end{equation*}
\end{definition}


\begin{definition}\label{D_Sublinear_Map}
A function $\varphi : X^{*} \to L_{1}(\mu)$ is said to be
\begin{itemize}
\item
\textit{positively homogeneous} 
if $\varphi(\alpha x^{*}) = \alpha \varphi(x)$ for every 
$\alpha \geq 0$ and $x^{*} \in X^{*}$,
\item
\textit{sub-additive} 
if $\varphi(x^{*} +y^{*}) \leq \varphi(x^{*}) + \varphi(y^{*})$ for all pairs 
$(x^{*}, y^{*}) \in X^{*} \times X^{*}$, where
\begin{equation*}
\begin{split}
f \leq g \Leftrightarrow f(t) \leq g(t)\text{ at almost all }t \in \Omega
\quad
(f,g \in L_{1}(\mu)),
\end{split}
\end{equation*}
\item
\textit{sublinear} 
if it is both positively homogeneous and sub-additive. 
\end{itemize}
\end{definition}


Thanks to James's theorem (\cite[Theorem 6]{JAMES}), 
the following one is proved in the same way as \cite[Proposition 2.4]{KAN}.

\begin{theorem}\label{T_CWK}
Let $M : \Sigma \to cwk(X)$ be a multimeasure. 
Then $M(\Sigma) = \bigcup_{E \in \Sigma} M(E)$ is a relatively weakly compact subset of $X$.
\end{theorem}


\section{The continuity of sublinear maps}

The continuity on $X^{*}$ of a 
sublinear map $\varphi: X^{*} \to L_{1}(\mu)$ ($\varphi: X^{*} \to \mathbb{R}$) 
is characterized by its continuity at $0 \in X^{*}$ 
with respect to 
a Hausdorff locally convex topology $\mathfrak{T}^{*}$ on $X^{*}$, 
see Lemma \ref{L_SublinearL1_Cont_Zero} (Lemma \ref{L_SublinearR_Cont_Zero}).
The notion of a bounded sublinear map is defined (Definition \ref{D_Bounded_Sublinear}). 
We show that the class of continuous sublinear maps 
coincides with 
the class of bounded sublinear maps  
(Lemma \ref{L_SublinearL1_Cont_Bounded} (Lemma \ref{L_SublinearR_Cont_Bounded})).


\begin{lemma}\label{L_SublinearL1_Cont_Zero}
Let $\mathfrak{T}^{*}$ be a Hausdorff locally convex topology in $X^{*}$ and 
let $\varphi : X^{*} \to L_{1}(\mu)$  be a sublinear function. 
Then the following statements are equivalent:
\begin{itemize}
\item[(i)]
$\varphi$ is continuous on $X^{*}$,
\item[(ii)]
$\varphi$ is continuous at $0 \in X^{*}$. 
\end{itemize}
\end{lemma}
\begin{proof}
It is clear that $(i) \Rightarrow (ii)$. 
To prove that  $(ii) \Rightarrow (i)$, 
we assume that $\varphi$ is continuous at $0 \in X^{*}$ and let $x_{0}^{*} \neq 0$.  
Then, given $\varepsilon >0$ there exists an absolutely convex $0$-neighborhood $U_{\varepsilon}$ in $X^{*}$ such that 
\begin{equation}\label{eqL_Sublinear_Cont_Zero.1}
\begin{split}
x^{*} \in U_{\varepsilon} \Rightarrow ||\varphi(x^{*})||_{1} \leq \varepsilon.
\end{split}
\end{equation}

Since for any $x^{*}, y^{*} \in X^{*}$, 
\begin{equation*}
\begin{split}
\varphi(x^{*}) \leq \varphi(x^{*}-y^{*}) + \varphi(y^{*})
\Rightarrow& 
\varphi(x^{*}) -\varphi(y^{*}) \leq \varphi(x^{*}-y^{*})
\\
\varphi(y^{*}) \leq \varphi(y^{*}-x^{*}) + \varphi(x^{*})
\Rightarrow& 
\varphi(y^{*}) -\varphi(x^{*}) \leq \varphi(y^{*}-x^{*}),
\end{split}
\end{equation*}
it follows that 
\begin{equation*}
\begin{split} 
-\varphi(y^{*}-x^{*})
\leq
\varphi(x^{*})-\varphi(y^{*})
\leq 
\varphi(x^{*}-y^{*}).
\end{split}
\end{equation*}
Thus,
\begin{equation*}
\begin{split}
-\varphi(y^{*}-x^{*})(t)
\leq
\varphi(x^{*})(t)-\varphi(y^{*})(t)
\leq 
\varphi(x^{*}-y^{*})(t)~\text{at almost all }t \in \Omega
\end{split}
\end{equation*}
and, consequently, 
\begin{equation*}
\begin{split}
|\varphi(x^{*})(t) -\varphi(y^{*})(t)|
\leq&
\max \{ |\varphi(y^{*}-x^{*})(t)|, |\varphi(x^{*}-y^{*})(t)|\}
\\
\leq&
|\varphi(y^{*}-x^{*})(t)|+ |\varphi(x^{*}-y^{*})(t)|
~\text{ at almost all }t \in \Omega.
\end{split}
\end{equation*}
Hence, 
\begin{equation}\label{eqL_Sublinear_Cont_Zero.2}
\begin{split}
||\varphi(x^{*}) -\varphi(y^{*})||_{1}
=&
\int_{\Omega} |\varphi(x^{*})(t) -\varphi(y^{*})(t)| d \mu
\\
\leq&
\int_{\Omega} |\varphi(y^{*}-x^{*})(t)| d \mu + \int_{\Omega} |\varphi(x^{*}-y^{*})(t)| d \mu 
\\
=&
||\varphi(y^{*}-x^{*})||_{1} + ||\varphi(x^{*}-y^{*})||_{1}.
\end{split}
\end{equation}
By \eqref{eqL_Sublinear_Cont_Zero.1} we have also
\begin{equation*}
\begin{split}
x^{*}-x^{*}_{0} \in U_{\varepsilon} 
\Rightarrow& 
x^{*}-x^{*}_{0} \in U_{\varepsilon}  \text{ and }x^{*}_{0}-x^{*} \in -U_{\varepsilon} = U_{\varepsilon}~(U_{\varepsilon}\text{ is circled})
\\
\Rightarrow& 
||\varphi(x^{*}-x^{*}_{0}||_{1} \leq \varepsilon \text{ and } ||\varphi(x^{*}_{0}-x^{*})||_{1} \leq \varepsilon.
\end{split}
\end{equation*}
The last result together with \eqref{eqL_Sublinear_Cont_Zero.2} yields
\begin{equation*}
\begin{split}
x^{*}-x^{*}_{0} \in U_{\varepsilon}  
\Rightarrow 
||\varphi(x^{*})-\varphi(x^{*}_{0})||_{1} \leq 2\varepsilon.
\end{split}
\end{equation*}
This means that $\varphi$ is continuous at $x^{*}_{0}$ and the proof is complete. 
\end{proof}


The following lemma can be proved in similar manner as Lemma \ref{L_SublinearL1_Cont_Zero}.

\begin{lemma}\label{L_SublinearR_Cont_Zero}
Let $\mathfrak{T}^{*}$ be a Hausdorff locally convex topology in $X^{*}$ and 
let $\varphi: X^{*} \to \mathbb{R}$ be a sublinear function. 
Then the following statements are equivalent:
\begin{itemize}
\item[(i)]
$\varphi$ is continuous on $X^{*}$,
\item[(ii)]
$\varphi$ is continuous at $0 \in X^{*}$. 
\end{itemize}
\end{lemma}


\begin{definition}\label{D_Bounded_Sublinear}
Let $\mathfrak{T}^{*}$ be a Hausdorff locally convex topology in $X^{*}$.
We say that a sublinear function $\varphi: X^{*} \to L_{1}(\mu)$ 
is \textit{bounded} if 
there exist $K>0$ and a continuous seminorm $p$ on $(X^{*}, \mathfrak{T}^{*})$ such that 
\begin{equation*}
\begin{split}
||\varphi(x^{*})||_{1} \leq K p(x^{*})
~\text{ for every }x^{*} \in X^{*}
\qquad
\left (
||f||_{1} = \int_{\Omega} |f(t)| d\mu, ~f\in L_{1}(\mu)
\right).
\end{split}
\end{equation*}
A sublinear function $\varphi: X^{*} \to \mathbb{R}$ is said to be \textit{bounded} if 
there exist $K>0$ and a continuous seminorm $p$ on $(X^{*}, \mathfrak{T}^{*})$ such that  
$|\varphi(x^{*})| \leq K p(x^{*})$ for all $x^{*} \in X^{*}$. 
\end{definition}


\begin{lemma}\label{L_SublinearL1_Cont_Bounded}
Let $\mathfrak{T}^{*}$ be a Hausdorff locally convex topology in $X^{*}$ 
and let $\varphi: X^{*} \to L_{1}(\mu)$ be a sublinear function. 
Then the following statements are equivalent:
\begin{itemize}
\item[(i)]
$\varphi$ is continuous on $X^{*}$, 
\item[(ii)]
$\varphi$ is a bounded sublinear function. 
\end{itemize}
\end{lemma}
\begin{proof}
$(i) \Rightarrow (ii)$   
Assume that $\varphi$ is continuous and let $\varepsilon =1$. 
Then there exists a continuous seminorm $p$ on $X^{*}$ such that 
for each $x^{*} \in X^{*}$, we have
\begin{equation}\label{eqL_SublinearL1_Cont_Bounded.1}
\begin{split}
p(x^{*}) \leq 1 \Rightarrow ||\varphi(x^{*})||_{1} \leq 1.
\end{split}
\end{equation}
Note that if $x^{*} \in X^{*}$ is such that $p(x^{*}) \neq 0$, then 
\begin{equation*}
\begin{split}
p \left ( \frac{1}{2 p(x^{*})} x^{*} \right ) = \frac{1}{2} < 1 
\Rightarrow 
 \left \| \frac{1}{2p(x^{*})}\varphi(x^{*}) \right \|_{1} \leq 1 
\Rightarrow 
||\varphi(x^{*})||_{1} \leq 2. p(x^{*}).
\end{split}
\end{equation*} 
Consequently, $||\varphi(x^{*})||_{1} \leq 2. p(x^{*})$ whenever $x^{*} \in X^{*}$ with $p(z^{*}) \neq 0$, 
and since for each $\varepsilon >0$, we have
\begin{equation*}
\begin{split}
p(x^{*}) \leq \varepsilon \Rightarrow 
||\varphi(x^{*})||_{1} \leq \varepsilon
\quad
(\text{see \eqref{eqL_SublinearL1_Cont_Bounded.1}}), 
\end{split}
\end{equation*}
it follows that $||\varphi(x^{*})||_{1} \leq 2. p(x^{*})$ for all $x^{*} \in X^{*}$. 
This means that $\varphi$ is a bounded sublinear function.

$(ii) \Rightarrow (i)$
Assume that $\varphi$ is a bounded sublinear function and let $\varepsilon>0$.
Then, there are a continuous seminorm $p$ on $X^{*}$ and $K >0$ such that 
$||\varphi(x^{*})||_{1} \leq K . p(x^{*})$ for all $x^{*} \in X^{*}$. 
Hence, 
\begin{equation*}
\begin{split}
x^{*} \in \frac{\varepsilon}{2K} U_{p}
\Rightarrow 
p(x^{*}) \leq \frac{\varepsilon}{2K}  
\Rightarrow 
||\varphi(x^{*})||_{1} \leq K . p(x^{*}) \leq K .\frac{\varepsilon}{2K} = \frac{\varepsilon}{2} < \varepsilon,
\end{split}
\end{equation*}
where $U_{p} = \{x^{*} \in X^{*}: p(x^{*}) \leq 1\}$. 
This means that $\varphi$ is continuous at $0 \in X^{*}$. 
Hence, by Lemma \ref{L_SublinearL1_Cont_Zero} the sublinear function $\varphi$ is continuous on $X^{*}$  and the proof is complete.
\end{proof}


The following lemma can be proved in similar way as Lemma \ref{L_SublinearL1_Cont_Bounded}.

\begin{lemma}\label{L_SublinearR_Cont_Bounded}
Let $\mathfrak{T}^{*}$ be a Hausdorff locally convex topology in $X^{*}$ and 
let $\varphi: X^{*} \to \mathbb{R}$ be a sublinear function. 
Then the following statements are equivalent:
\begin{itemize}
\item[(i)]
$\varphi$ is continuous on $X^{*}$, 
\item[(ii)]
$\varphi$ is a bounded sublinear function. 
\end{itemize}
\end{lemma}


\section{Characterizations of Pettis integrability in locally convex spaces}

In this section, the main results are Lemma \ref{L_Pettis_BCC} and  Theorem \ref{T_Pettis_Chr_CWK}. 
Let us start with an auxiliary lemma.


\begin{lemma}\label{L_Convex_Hull_W_Compact}
Let $(X,\mathfrak{T})$ be a complete Hausdorff locally convex topological vector space 
and let $A \in \mathcal{P}(X)$ be a weakly compact set. 
Then the closed convex hull $\overline{\text{co}}(A)$ of $A$ is weakly compact. 
\end{lemma}
\begin{proof}
Since $C = \overline{\text{co}}(A)$ is a closed subset of the complete Hausdorff locally convex space $X$, 
$C$ is $\mathfrak{T}$-complete. 
Thanks to \cite[Theorem 11.4, p.189]{SCH}, 
it is enough to prove that $C$ is  complete for the Mackey topology $\tau(X,X^{*})$. 
To do this, 
let $(x_{\alpha})_{\alpha \in \mathcal{A}}$ be a $\tau(X,X^{*})$-Cauchy net with terms in $C$.

Since Mackey topology $\tau(X,X^{*})$ is the topology of uniform convergence on the family 
$\mathcal{M}(X^{*},X)$ of all absolutely convex and $\sigma(X^{*},X)$-compact subsets of $X^{*}$,  
it is clear that the family 
$\{H^{o}: H \in \mathcal{M}(X^{*},X)\}$ 
is a $0$-neighborhood base under the Mackey topology $\tau(X, X^{*})$, 
where 
$H^{o} = \left \{x \in X: \sup_{x^{*} \in H} |x^{*}(x)| \leq 1 \right \}$. 
Then, 
given $H \in \mathcal{M}(X^{*},X)\}$ 
there exists $\alpha_{H} \in \mathcal{A}$ such that 
\begin{equation}\label{eqL_Convex_Hull_W_Compact.1}
\begin{split}
\alpha, \beta \geq \alpha_{H} 
\Rightarrow 
x_{\alpha} - x_{\beta} \in H^{o}
\Rightarrow 
\sup_{x^{*} \in H} |x^{*}(x_{\alpha}) - x^{*}(x_{\beta})| \leq 1
\end{split}
\end{equation}

From the fact that $\mathfrak{T} \leq \tau(X, X^{*})$, 
it follows that 
$(x_{\alpha})_{\alpha \in \mathcal{A}}$ is a $\mathfrak{T}$-Cauchy net with terms in $C$. 
Hence, the net $(x_{\alpha})_{\alpha \in \mathcal{A}}$ converges to $x_{0} \in C$ 
under topology $\mathfrak{T}$, 
and since $\sigma(X,X^{*}) \leq \mathfrak{T}$, 
we have
\begin{equation}\label{eqL_Convex_Hull_W_Compact.2}
\begin{split}
\lim_{\beta \in \mathcal{A}} x^{*}(x_{\beta}) = x^{*}(x_{0})
~\text{ for every }x^{*} \in X^{*}. 
\end{split}
\end{equation}

Fix an arbitrary $\alpha \geq \alpha_{H}$. 
Then, by \eqref{eqL_Convex_Hull_W_Compact.2} and \eqref{eqL_Convex_Hull_W_Compact.1} 
we get
\begin{equation*}
\begin{split}
\lim_{\beta \in \mathcal{A}} |x^{*}(x_{\alpha}) - x^{*}(x_{\beta})| 
= 
|x^{*}(x_{\alpha}) - x^{*}(x_{0})| \leq 1
~\text{ for all }x^{*} \in H.
\end{split}
\end{equation*}
Since $\alpha$ was arbitrary, the last result yields 
\begin{equation*}
\begin{split}
\alpha \geq \alpha_{H} 
\Rightarrow 
\sup_{x^{*} \in H} |x^{*}(x_{\alpha}) - x^{*}(x_{0})| \leq 1
\Rightarrow 
x_{\alpha} - x_{0} \in H^{o}.
\end{split}
\end{equation*}
Thus, the $\tau(X,X^{*})$-Cauchy net
$(x_{\alpha})_{\alpha \in \mathcal{A}}$ with terms in $C$ 
converges to $x_{0} \in C$ under $\tau(X,X^{*})$. 
This means that $C$ is  $\tau(X,X^{*})$-complete, 
and this completes the proof.
\end{proof}

Since by \cite[5.2, p.26]{SCH} the circled hull of every weakly compact subset of $X$ is
weakly compact, 
Lemma \ref{L_Convex_Hull_W_Compact} continues to hold with $C$ the closed absolutely convex 
hull of $A$, 
i.e., $C = \overline{\text{aco}}(A)$.


\begin{definition}\label{D_V_Gamma}
Given a multifunction $\Gamma : \Omega \to cc(X)$ scalarly integrable 
we can define two sublinear maps $\varphi_{\Gamma} : X^{*} \to L_{1}(\mu)$ and 
$\varphi_{\Gamma}^{A}: X^{*} \to \mathbb{R}$ for every $A \in \Sigma$ as follows
\begin{equation*}
\begin{split}
\varphi_{\Gamma}(x^{*}) = \sigma(x^{*}, \Gamma(\cdot)), 
\quad
\varphi_{\Gamma}^{A}(x^{*}) = \int_{A} \sigma(x^{*}, \Gamma(t)) d \mu.
\end{split}
\end{equation*}
\end{definition}


The next lemma is easily proved by \cite[Theorem II-16, p.48]{CAST}.

\begin{lemma}\label{L_Pettis_Chr_BCC}
Let $\Gamma : \Omega \to cc(X)$ be a scalarly integrable multifunction.
Then the following statements are equivalent:
\begin{itemize}
\item[(i)]
$\Gamma$ is Pettis integrable in $bcc(X)$, 
\item[(ii)]
for each $A \in \Sigma$ the sublinear map 
$\varphi^{A}_{\Gamma}$ is $\sigma(X^{*},X)$-lower semicontinuous.
\end{itemize}
\end{lemma}
\begin{proof}
$(i) \Rightarrow (ii)$ 
Assume that $\Gamma$ is Pettis integrable in $bcc(X)$ 
with $M(A) = (P)\int_{A} \Gamma(t) d \mu$ for all $A \in \Sigma$. 
Then, 
\begin{equation*}
\begin{split}
\varphi^{A}_{\Gamma}(x^{*}) = \sigma(x^{*}, M(A)) = \int_{A} \sigma(x^{*}, \Gamma(t)) d \mu
\qquad
(A \in \Sigma, x^{*} \in X^{*}). 
\end{split}
\end{equation*} 
Therefore, by \cite[Theorem II-16, p.48]{CAST} 
the function $\varphi^{A}_{\Gamma}$ is $\sigma(X^{*},X)$-lower semi-continuous.

$(ii) \Rightarrow (i)$ 
Assume that for each $A \in \Sigma$ the sublinear function 
$\varphi^{A}_{\Gamma}$ is $\sigma(X^{*},X)$-lower semi-continuous. 
Then,  by \cite[Theorem II-16, p.48]{CAST} for each $A \in \Sigma$ there exists a unique 
$M(A) \in cc(X)$ such that 
\begin{equation*}
\begin{split}
\sigma(x^{*}, M(A))
= 
\varphi^{A}_{\Gamma}(x^{*}) 
=
\int_{A} \sigma(x^{*}, \Gamma(t)) d \mu  \in \mathbb{R}
~\text{ for all }x^{*} \in X^{*}.
\end{split}
\end{equation*}
Then, for each $A \in \Sigma$, the set $M(A)$ is  a $\sigma(X,X^{*})$-bounded subset of $X$. 
Hence, by \cite[Corollary 2, p.132]{SCH} (or, \cite[Theorem 1, p.67]{ROB}) the set $M(A)$ is a bounded subset of $X$. 
Consequently, $\Gamma$ is Pettis integrable in $bcc(X)$ and the proof is complete.
\end{proof}

The following example shows that there is 
a multifunction $\Gamma : \Omega \to bcc(X)$ Pettis integrable in $bcc(X)$ such that 
$\varphi_{\Gamma}$ is  $\beta(X^{*},X)$-$||\cdot||_{1}$-continuous,  
but it  not $\tau(X^{*},X)$-$||\cdot||_{1}$-continuous.

\begin{example}\label{EXE_BCC}
Let $X$ be a Banach space and let $f: \Omega \to X$ be a Pettis integrable function 
with $\nu(A) = (P)\int_{A} f(t) d \mu$  for all $A \in \Sigma$. 
According to \cite[Example 1.11]{MUS_Pettis_MultiF} the multifunction $\Gamma : \Omega \to bcc(X)$ defined by
\begin{equation*}
\begin{split}
\Gamma(t) = f(t) + r(t) B_{X}, 
\quad
(B_{X} = \{x \in X: ||x|| \leq 1\}\text{ and }
r: \Omega \to (0,+\infty)\text{ is Lebesgue integrable}), 
\end{split}
\end{equation*}
is Pettis integrable in $bcc(X)$ with 
\begin{equation*}
\begin{split}
M(A) = (P) \int_{A} \Gamma(t) d \mu = (P)\int_{A} f(t) d \mu + 
\left ( \int_{A} r(t) d \mu \right ) B_{X}. 
\end{split}
\end{equation*}
Hence, 
\begin{equation*}
\begin{split}
|\sigma(x^{*}, M(A))| \leq ||x^{*}|| \left ( ||\nu(A)|| + \int_{A} r(t) d \mu \right )
\qquad
(x^{*} \in X^{*}, A\in \Sigma), 
\end{split}
\end{equation*}
and by \cite[Lemma III.1.5, p.97]{DUN} we have also
\begin{equation*}
\begin{split}
||\varphi_{\Gamma}(x^{*})||_{1} = \int_{\Omega} \sigma(x^{*},\Gamma(t)) d \mu \leq 4 \sup_{A \in \Sigma} |\sigma(x^{*}, M(A))|
\quad
(x^{*} \in X^{*}). 
\end{split}
\end{equation*}
Since $\nu(\Sigma) =\{\nu(A):A \in \Sigma\}$ 
is relatively weakly compact 
(see \cite[Corollary 7, p.14]{DIES}) 
and  $0 \leq \int_{A} r(t) d \mu \leq \int_{\Omega} r(t) d \mu$, 
it follows that 
\begin{equation*}
\begin{split}
||\varphi_{\Gamma}(x^{*})||_{1} \leq 4 ||x^{*}|| \left ( \sup_{A \in \Sigma} ||\nu(A)|| + \int_{A} r(t) d \mu \right ).
\end{split}
\end{equation*}
Since the strong topology $\beta(X^{*},X)$ coincides with the topology generated by the operator norm, 
the last result yields that $\varphi_{\Gamma}$ is  $\beta(X^{*},X)$-$||\cdot||_{1}$-continuous.

Moreover, by \cite[Example 1.11]{MUS_Pettis_MultiF} if $X = \ell_{1}$, 
then $\Gamma$ is Pettis integrable in $bcc(\ell_{1})$ 
but $\varphi_{\Gamma}$ is not $\tau(\ell_{1}^{*},\ell_{1})$-$\sigma(L_{1},L_{1}^{*})$-continuous 
on $B(\ell_{\infty})$, ($\ell_{\infty} \equiv \ell_{1}^{*}$). 
Hence, $\varphi_{\Gamma}$ is not $\tau(\ell_{1}^{*},\ell_{1})$-$\sigma(L_{1},L_{1}^{*})$-continuous 
and, consequently, it is not $\tau(\ell_{1}^{*},\ell_{1})$-$||\cdot||_{1}$-continuous.
\end{example}

The following lemma shows that for each multifunction $\Gamma : \Omega \to cc(X)$ Pettis integrable in $bcc(X)$,  
the map $\varphi_{\Gamma}$ is  $\beta(X^{*},X)$-$||\cdot||_{1}$-continuous.

\begin{lemma}\label{L_Pettis_BCC}
If a multifunction $\Gamma : \Omega \to cc(X)$ is Pettis integrable in $bcc(X)$ 
with 
$M(E) = (P)\int_{E} \Gamma(t) d \mu$ for all $E \in \Sigma$, 
then
\begin{itemize}
\item[(i)]
$M(\Sigma) = \bigcup_{E \in \Sigma} M(E)$ is bounded subset of $X$, 
\item[(ii)]
$\varphi_{\Gamma}$ is $\beta(X^{*},X)$-$||\cdot||_{1}$-continuous.
\end{itemize}
\end{lemma}
\begin{proof}
$(i)$ 
For each $x^{*} \in X^{*}$ the real valued countable additive measure $\nu_{x^{*}} : \Sigma \to \mathbb{R}$ define by 
\begin{equation*}
\begin{split}
\nu_{x^{*}}(E) = \sigma(x^{*}, M(E)) = \int_{E} \sigma(x^{*}, \Gamma(t)) d \mu
~\text{ for every }E \in \Sigma 
\end{split}
\end{equation*}
is bounded (see \cite[Lemma III.4.4, p.141]{DUN}). 
Then, for each $x^{*} \in X^{*}$ 
there exists $K_{x^{*}} >0$ such that 
$|\nu_{x^{*}}(E)|\leq K_{x^{*}}$ for all $E \in \Sigma$, 
and since
\begin{equation*}
\begin{split}
\inf x^{*} C = - \sup (-x^{*}C)
\quad 
(C \in cc(X)),
\end{split}
\end{equation*}
it follows that 
\begin{equation*} 
\begin{split}
|\inf x^{*} M(E)| = |\sup((-x^{*})M(E))|
=
|\sigma(-x^{*}, M(E))| \leq  K_{-x^{*}} 
~\text{ for all }E \in \Sigma~(x^{*} \in X^{*}).
\end{split}
\end{equation*}
Hence, 
\begin{equation*}
\begin{split}
-K_{-x^{*}} \leq \inf x^{*}M(E) \leq \sup x^{*}M(E) \leq K_{x^{*}}
~\text{ for all } E \in \Sigma~(x^{*} \in X^{*}).
\end{split}
\end{equation*}
Consequently, 
\begin{equation*} 
\begin{split}
-K_{-x^{*}} \leq \inf_{E \in \Sigma} [\inf x^{*}M(E)] \leq \sup_{E \in \Sigma} [\sup x^{*}M(E)] \leq K_{x^{*}}
~\text{ for every } x^{*} \in X^{*}.
\end{split}
\end{equation*}
The last result together with 
\begin{equation*}
\begin{split}
\inf_{E \in \Sigma} [\inf x^{*}M(E)] = \inf x^{*}M(\Sigma),
\quad
\sup_{E \in \Sigma} [\sup x^{*}M(E)] = \sup x^{*}M(\Sigma)
\end{split}
\end{equation*}
yields 
\begin{equation*}
\begin{split}
-K_{-x^{*}} \leq \inf x^{*}M(\Sigma) \leq \sup x^{*}M(\Sigma) \leq K_{x^{*}} 
~\text{ for every }x^{*} \in X^{*}.  
\end{split}
\end{equation*}
This means that $M(\Sigma)$ is a $\sigma(X,X^{*})$-bounded subset of $X$. 
Hence, by \cite[Corollary 2, p.132]{SCH} (or, \cite[Theorem 1, p.67]{ROB}) the set $M(\Sigma)$ is a bounded subset of $X$.

$(ii)$    
Since each real valued measure $\nu_{x^{*}}$ is bounded, 
we obtain by \cite[Lemma III.1.5, p.97]{DUN} that $\nu_{x^{*}}$ is of bounded variation and 
$|\nu_{x^{*}}|(\Omega) \leq 4\sup_{E \in \Sigma} |\sigma(x^{*}, M(E))|$, and since 
\begin{equation*}
\begin{split}
|\nu_{x^{*}}|(\Omega) = \int_{\Omega} |\sigma(x^{*}, \Gamma(t))| d \mu =||\varphi_{\Gamma}(x^{*})||_{1},
\end{split}
\end{equation*} 
it follows that 
\begin{equation*} 
\begin{split}
|\nu_{x^{*}}|(\Omega) = \int_{\Omega} |\sigma(x^{*}, \Gamma(t))| d \mu =||\varphi_{\Gamma}(x^{*})||_{1} 
\leq  4\sup_{E \in \Sigma} |\sigma(x^{*}, M(E))|
\quad
(x^{*} \in X^{*}).
\end{split}
\end{equation*}
Hence, by inequalities 
\begin{equation*}  
\begin{split}
\sup_{E \in \Sigma} |\sigma(x^{*}, M(E))| 
\leq
\sup_{E \in \Sigma} \sup_{x \in M(E)} |x^{*}(x)|
\leq
\sup_{x \in M(\Sigma)} |x^{*}(x)|
\quad
(x^{*} \in X^{*}),
\end{split}
\end{equation*}
we obtain
\begin{equation}\label{eqT_Pettis_Chr_BCC.KEY}
\begin{split}
||\varphi_{\Gamma}(x^{*})||_{1} \leq 4\sup_{x \in M(\Sigma)} |x^{*}(x)| = 4 p_{M(\Sigma)}(x^{*})
~\text{ for all }x^{*} \in X^{*}, 
\end{split}
\end{equation}
where $p_{M(\Sigma)}$ is the Minkowski functional of 
the absolute polar $M(\Sigma)^{o}$ of $M(\Sigma)$. 
Thus, $p_{M(\Sigma)}$ is a $\beta(X^{*},X)$-continuous seminorm on $X^{*}$  
and by Definition \ref{D_Bounded_Sublinear} follows that $\varphi_{\Gamma}$ is a bounded sublinear function. 
Therefore, by Lemma \ref{L_SublinearL1_Cont_Bounded} the sublinear function $\varphi_{\Gamma}$ 
is $\beta(X^{*},X)$-$||\cdot||_{1}$ continuous, and this completes the proof.  
\end{proof}


\begin{theorem}\label{T_Pettis_Chr_CWK}
Let $\Gamma : \Omega \to cc(X)$ be a scalarly integrable multifunction. 
Then the following statements are equivalent:
\begin{itemize}
\item[(i)]
$\Gamma$ is Pettis integrable in $cwk(X)$, 
\item[(ii)]
$\varphi_{\Gamma}$ is $\tau(X^{*},X)$-$||\cdot||_{1}$ continuous,
\item[(iii)]
$\varphi_{\Gamma}$ is $\tau(X^{*},X)$-$\sigma(L_{1},L_{1}^{*})$-continuous,
\item[(iv)]
for each $A \in \Sigma$ the map $\varphi_{\Gamma}^{A}$ is  $\tau(X^{*},X)$-continuous.
\end{itemize}
\end{theorem}
\begin{proof}
$(i) \Rightarrow (ii)$    
Assume that $\Gamma$ is Pettis integrable in $cwk(X)$ 
with $M(A) = (P)\int_{A} \Gamma(t) d \mu$ for all $A \in \Sigma$. 
Then, Theorem \ref{T_CWK} together with Lemma \ref{L_Convex_Hull_W_Compact} yields that 
the $\sigma(X,X^{*})$-closed absolutely convex hull 
$C = \overline{\text{aco}(M(\Sigma))}^{\sigma(X,X^{*})}$ 
of $M(\Sigma)$ is a weakly compact subset of $X$.  
By \eqref{eqT_Pettis_Chr_BCC.KEY} in Lemma \ref{L_Pettis_BCC} we have also
\begin{equation*} 
\begin{split}
||\varphi_{\Gamma}(x^{*})||_{1} 
\leq  4 p_{M(\Sigma)}(x^{*})
\leq 4 p_{C}(x^{*})
~\text{ for all }x^{*} \in X^{*}.
\end{split}
\end{equation*}
This means that the sublinear function $\varphi_{\Gamma}$ is bounded and hence 
by Lemma \ref{L_SublinearL1_Cont_Bounded} the map $\varphi_{\Gamma}$ 
is $\tau(X^{*},X)$-$||\cdot||_{1}$ continuous.

Since the weak topology $\sigma(L_{1},L_{1}^{*})$ in $L_{1}(\mu)$  
is coarser then the topology generated by norm $||\cdot||_{1}$, 
it follows that $(ii) \Rightarrow (iii)$.

$(iii) \Rightarrow (iv)$ 
Assume that $\varphi_{\Gamma}$ is $\tau(X^{*},X)$-$\sigma(L_{1},L_{1}^{*})$ continuous 
and let $\kappa: L_{\infty}(\mu) \to (L_{1}(\mu))^{*}$ be the isomorphism isometric between $(L_{1}(\mu))^{*}$ and $L_{\infty}(\mu)$ defined by 
\cite[Theorem IV.8.5]{DUN}. 
Then, 
\begin{equation*}
\begin{split}
\langle \varphi_{\Gamma}(x^{*}), \kappa(g) \rangle = \int_{\Omega} \sigma(x^{*}, \Gamma(t)) g(t) d \mu 
\quad
(x^{*} \in X^{*}, g \in L_{\infty}(\mu))
\end{split}
\end{equation*}
and, in particular, 
\begin{equation*}
\begin{split}
\langle \varphi_{\Gamma}(x^{*}), \kappa(\mathbbm{1}_{A}) \rangle = \int_{A} \sigma(x^{*}, \Gamma(t))d \mu =\varphi_{\Gamma}^{A}(x^{*})
\quad
(x^{*} \in X^{*}, A \in \Sigma).
\end{split}
\end{equation*}
Consequently, $\varphi_{\Gamma}^{A} = \kappa(\mathbbm{1}_{A}) \circ \varphi_{\Gamma}$ 
is $\tau(X^{*},X)$-continuous.

$(iv) \Rightarrow (i)$ 
Assume that $(iv)$ holds. 
Then, for each $r \in \mathbb{R}$ for each $A \in \Sigma$ 
the convex set 
$H(r,A) = \{x^{*} \in X^{*}: \psi^{A}_{\Gamma}(x^{*}) \leq r\}$ is $\tau(X^{*},X)$-closed, 
and since
\begin{equation*}
\begin{split}
H(r,A) 
= \overline{H(r,A)}^{ \tau(X^{*},X)}
= \overline{H(r,A)}^{ \sigma(X^{*},X)}
\qquad
(\text{cf. \cite[Proposition 8, p.34]{ROB}},
\end{split}
\end{equation*} 
it follows that $H(r,A)$ is $\sigma(X^{*},X)$-closed. 
This means that $\psi^{A}_{\Gamma}$ is $\sigma(X^{*},X)$-lower semicontinuous, 
cf. \cite[\S 6, 2(3), p.40]{KOTHE}. 
Consequently, by \cite[Theorem II-16, p.48]{CAST} there exists a unique $M(A) \in cc(X)$ such that 
\begin{equation*}
\begin{split}
\psi^{A}_{\Gamma}(x^{*}) = \sigma(x^{*}, M(A)) = \int_{A} \sigma(x^{*},\Gamma(t)) d \mu
\qquad
(x^{*} \in X^{*}, A \in \Sigma).
\end{split}
\end{equation*}
It remains to show that $M(A) \in cwk(X)$ for all $A \in \Sigma$.

Note that for a fixed $A \in \Sigma$, 
the sublinear map $\psi^{A}_{\Gamma}(\cdot) = \sigma(\cdot, M(A))$ is 
$\tau(X^{*},X)$-continuous. 
Then, given $\varepsilon =1$ there exists $H \in \mathcal{M}(X,X^{*})$ such that 
\begin{equation}\label{eqPolarMackey}
\begin{split}
x^{*} \in H^{o} 
= 
\left \{ z^{*} \in X^{*} : \sup_{x \in H} |z^{*}(x)| \leq 1 \right \} 
\Rightarrow |\sigma(x^{*}, M(A))| \leq 1 \text{ and }|\sigma(-x^{*}, M(A))| \leq 1,
\end{split}
\end{equation}
because 
$H^{o}$ is a circled set. 
Since
\begin{equation*}
\begin{split}
- \sigma(-x^{*}, M(A)) =\inf x^{*} M(A) \leq x^{*}(x) \leq \sup x^{*} M(A) =\sigma(x^{*}, M(A))
~\text{ for every } x \in M(A),
\end{split}
\end{equation*}
we obtain  
\begin{equation*}
\begin{split}
\sup_{x\in M(A)} |x^{*}(x)| \leq \max \{ |\sigma(-x^{*}, M(A))|, |\sigma(x^{*}, M(A))|\} 
\end{split}
\end{equation*}
The last result together with \eqref{eqPolarMackey} yields that 
\begin{equation*}
\begin{split}
x^{*} \in H^{o} \Rightarrow \sup_{x\in M(A)} |x^{*}(x)| \leq 1 \Rightarrow x^{*} \in 
M(A)^{o}= \left \{z^{*} \in X^{*} : \sup_{x \in M(A)} |z^{*}(x)| \leq 1 \right \}.
\end{split}
\end{equation*}
It follows by Bipolar Theorem (cf. \cite[Theorem 1.5, p.126]{SCH}) that
\begin{equation*}
\begin{split}
H^{o} \subset M(A)^{o} 
\Rightarrow 
M(A)^{oo} \subset H^{oo} 
\Rightarrow
\overline{\text{co}}^{\sigma(X,X^{*})}(M(A) \cup \{\theta\}) \subset H
\Rightarrow
M(A) = \overline{M(A)}^{\sigma(X,X^{*})} \subset H.  
\end{split}
\end{equation*}
Consequently, $M(A)$ is $\sigma(X,X^{*})$-compact. 
Thus, $\Gamma$ is Pettis integrable in $cwk(X)$, and this completes the proof. 
\end{proof}

Theorem \ref{T_Pettis_Chr_CWK} is the key to prove the following theorem.

\begin{theorem}\label{T_Convergence_Theorem}
Let $\Gamma: \Omega \to cc(X)$ be a multifunction 
and let $\{\Gamma_{n} : \Omega \to cc(X): n \in \mathbb{N}\}$ be a sequence of scalarly integrable multifunctions. 
Assume that 
\begin{itemize}
\item
there exists $U \in \mathcal{M}(X,X^{*})$ such that the  family $\{\sigma(x^{*}, \Gamma_{n}(\cdot)): x^{*} \in U^{o}, n \in \mathbb{N}\}$ 
is uniformly integrable, 
\item
for each $x^{*} \in X^{*}$, 
the sequence $(\sigma(x^{*}, \Gamma_{n}))$ converges in measure to $\sigma(x^{*}, \Gamma)$. 
\end{itemize}
Then 
\begin{itemize}
\item[(i)]
for each $n \in \mathbb{N}$ the multifunction $\Gamma_{n}$ is Pettis integrable in $cwk(X)$,
\item[(ii)]
$\Gamma$ is Pettis integrable in $cwk(X)$ and 
\begin{equation*}
\begin{split}
\text{(CT)}
\qquad
\lim_{n \to \infty} 
\sigma \left ( x^{*}, (P) \int_{A} \Gamma_{n}(t) d \mu \right )
=
\sigma \left ( x^{*}, (P) \int_{A} \Gamma(t) d \mu \right )
\quad
(x^{*} \in X^{*}, A \in \Sigma). 
\end{split}
\end{equation*}
\end{itemize}
\end{theorem}
\begin{proof}
$(i)$ 
Let $\varphi_{\Gamma_{n}}$ be the sublinear map defined by $\Gamma_{n}$ 
for each $n \in \mathbb{N}$, i.e., 
\begin{equation*}
\begin{split}
\varphi_{\Gamma_{n}} : X^{*} \to L_{1}(\mu), 
\quad
\varphi_{\Gamma_{n}}(x^{*}) = \sigma(x^{*}, \Gamma_{n}(\cdot)).
\end{split}
\end{equation*}
The Minkowski functional $p_{U}$ of $U^{o}$ is defined by 
\begin{equation*}
\begin{split}
p_{U} : X^{*} \to [0,+\infty), 
\quad
p_{U}(x^{*}) = \sup_{x \in U} |x^{*}(x)|. 
\end{split}
\end{equation*}
Since the absolute polar $U^{o}$ of $U \in \mathcal{M}(X,X^{*})$ is a $0$-neighborhood in $X^{*}$ 
under the Mackey topology $\tau(X^{*},X)$, 
$p_{U}$ is a continuous seminorm on $X^{*}$ with respect to $\tau(X^{*},X)$-topology.

From the fact that the family $\{\sigma(x^{*}, \Gamma_{n}(\cdot)): x^{*} \in U^{o}, n \in \mathbb{N}\}$ 
is uniformly integrable, 
it follows that there exists $K>0$ such that 
\begin{equation*}
\begin{split}
||\sigma(x^{*}, \Gamma_{n}(\cdot))||_{1} = ||\varphi_{\Gamma_{n}}(x^{*})||_{1} \leq K
~\text{ for all }x^{*} \in U^{o}\text{ and }n \in \mathbb{N},
\end{split}
\end{equation*}
and
\begin{equation*}
\begin{split}
\lim_{\mu(A) \to 0} \int_{A} |\sigma(x^{*}, \Gamma_{n}(t))| d \mu = 0
~\text{ uniformly for }x^{*} \in U^{o}\text{ and }n \in \mathbb{N}.
\end{split}
\end{equation*} 
Then, for each $x^{*} \in X^{*}$, we have
\begin{equation}\label{eqT_Convergence_Theorem.1}
\begin{split}
p_{U}(x^{*}) < 1 \Rightarrow x^{*} \in U^{o} \Rightarrow ||\varphi_{\Gamma_{n}}(x^{*})||_{1} \leq K\text{ for all }n \in \mathbb{N}.
\end{split}
\end{equation}
Hence, for each $x^{*} \in X^{*}$ such that $p_{U}(x^{*}) \neq 0$, 
\begin{equation*}
\begin{split}
||\varphi_{\Gamma_{n}}(x^{*})||_{1} \leq 2 K p_{U}(x^{*})\text{ for all }n \in \mathbb{N},
\end{split}
\end{equation*}
and since for each $\varepsilon >0$, 
\begin{equation*}
\begin{split}
p_{U}(x^{*}) < \varepsilon \Rightarrow ||\varphi_{\Gamma_{n}}(x^{*})||_{1} \leq \varepsilon K\text{ for all }n \in \mathbb{N}
\qquad
(\text{see \eqref{eqT_Convergence_Theorem.1}}), 
\end{split}
\end{equation*}
it follows that 
\begin{equation*}
\begin{split}
||\varphi_{\Gamma_{n}}(x^{*})||_{1} \leq 2 K p_{U}(x^{*})\text{ for all }x^{*} \in X^{*}\text{ and }n \in \mathbb{N}.
\end{split}
\end{equation*}
The last result together with Definition \ref{D_Bounded_Sublinear} 
and Lemma \ref{L_SublinearL1_Cont_Bounded} yields that 
$\varphi_{\Gamma_{n}}$ is $\tau(X^{*},X)$-$||\cdot||_{1}$ continuous 
for every $n \in \mathbb{N}$. 
Consequently, by Theorem \ref{T_Pettis_Chr_CWK} each multifunction 
$\Gamma_{n}$ is Pettis integrable in $cwk(X)$.

$(ii)$ 
Note that for each $x^{*} \in U^{o}$ we have
\begin{itemize}
\item
the family $\{ \sigma(x^{*},\Gamma_{n}) : n \in \mathbb{N} \}$ is uniformly integrable,
\item
the sequence $(\sigma(x^{*}, \Gamma_{n}))$ converges in measure to $\sigma(x^{*}, \Gamma)$.
\end{itemize}
Then by Vitali convergence theorem (cf. \cite[Theorem III.3.6, p.122]{DUN}) we obtain 
that $\Gamma$ is scalarly integrable and 
\begin{equation}\label{eqT_Convergence_Theorem.2}
\begin{split}
\lim_{n \to \infty} ||\varphi_{\Gamma_{n}}(x^{*}) - \varphi_{\Gamma}(x^{*})||_{1}
=
\lim_{n \to \infty} ||\sigma(x^{*},\Gamma_{n}(\cdot)) - \sigma(x^{*},\Gamma(\cdot))||_{1} = 0
~\text{ for every }x^{*} \in U^{o}.
\end{split}
\end{equation}
Then,  
\begin{equation*}
\begin{split}
\lim_{n \to \infty} ||\varphi_{\Gamma_{n}}(x^{*})||_{1} = ||\varphi_{\Gamma}(x^{*})||_{1}
~\text{ for every }x^{*} \in U^{o}
\end{split}
\end{equation*}
and hence 
$||\varphi_{\Gamma}(x^{*})||_{1} \leq K$ for all $x^{*} \in U^{o}$ (see \eqref{eqT_Convergence_Theorem.1}).
It follows that  
\begin{equation*}
\begin{split}
||\varphi_{\Gamma}(x^{*})||_{1} \leq 2 K p_{U}(x^{*})\text{ for all }x^{*} \in X^{*}.
\end{split}
\end{equation*}
Thus, $\varphi_{\Gamma}$ is $\tau(X^{*},X)$-$||\cdot||_{1}$ continuous and, consequently 
by Theorem \ref{T_Pettis_Chr_CWK}  the multifunction 
$\Gamma$ is Pettis integrable in $cwk(X)$.

Note that given $x^{*} \in U^{o}, n \in \mathbb{N}$ and $A \in \Sigma$, 
\begin{equation*}
\begin{split}
\left |\int_{A} \sigma(x^{*}, \Gamma_{n}(t)) d \mu 
- \int_{A} \sigma(x^{*}, \Gamma(t)) d \mu \right |
\leq&
\int_{A} |\sigma(x^{*}, \Gamma_{n}(t))-\sigma(x^{*}, \Gamma(t))| d \mu
\\
\leq&
\int_{\Omega} |\sigma(x^{*}, \Gamma_{n}(t))-\sigma(x^{*}, \Gamma(t))| d \mu
=||\varphi_{\Gamma_{n}}(x^{*}) - \varphi_{\Gamma}(x^{*})||_{1}.
\end{split}
\end{equation*}
It follows by \eqref{eqT_Convergence_Theorem.2} that
\begin{equation*}
\begin{split}
\lim_{n \to \infty} \int_{A} \sigma(x^{*}, \Gamma_{n}(t)) d \mu = \int_{A} \sigma(x^{*}, \Gamma(t)) d \mu
\qquad
(x^{*} \in U^{o}, A \in \Sigma). 
\end{split}
\end{equation*}
If we fix an arbitrary $z^{*} \not\in U^{o}$, 
then there exists $\lambda >0$ such that $\lambda x^{*} \in U^{o}$ and, therefore,
\begin{equation*}
\begin{split}
\lambda \lim_{n \to \infty} \int_{A} \sigma(z^{*}, \Gamma_{n}(t)) d \mu
=&
\lim_{n \to \infty} \int_{A} \sigma(\lambda z^{*}, \Gamma_{n}(t)) d \mu 
= 
\int_{A} \sigma(\lambda z^{*}, \Gamma(t)) d \mu
\\
=&
\lambda \int_{A} \sigma(z^{*}, \Gamma(t)) d \mu
\qquad
(z^{*} \in X^{*}, A \in \Sigma). 
\end{split}
\end{equation*}
This means that \text{(CT)} holds, and this completes the proof.
\end{proof}

The following corollary shows a version of the Lebesgue dominated convergence theorem 
for the scalarly integrable multifunctions.

\begin{corollary}\label{C_Dominated_Convergence_Theorem}
Let $\Gamma: \Omega \to cc(X)$ be a multifunction 
and let $\{\Gamma_{n} : \Omega \to cc(X): n \in \mathbb{N}\}$ be a sequence of scalarly integrable multifunctions. 
Assume that 
\begin{itemize}
\item
there exist $g \in L_{1}(\mu)$ and $U \in \mathcal{M}(X,X^{*})$ 
such that for each $x^{*} \in U^{o}$ and for each $n \in \mathbb{N}$, 
\begin{equation*}
\begin{split}
|\sigma(x^{*}, \Gamma_{n}(t))| \leq g(t)
~\text{ at almost all }t \in \Omega,
\end{split}
\end{equation*}
\item
for each $x^{*} \in X^{*}$, 
the sequence $(\sigma(x^{*}, \Gamma_{n}))$ converges in measure to $\sigma(x^{*}, \Gamma)$. 
\end{itemize}
Then 
\begin{itemize}
\item[(i)]
for each $n \in \mathbb{N}$ the multifunction $\Gamma_{n}$ is Pettis integrable in $cwk(X)$,
\item[(ii)]
$\Gamma$ is Pettis integrable in $cwk(X)$ and \text{(CT)} holds.
\end{itemize}
\end{corollary}
\begin{proof}
Thanks to Theorem \ref{T_Convergence_Theorem}, 
it is enough to prove that the family 
$\mathcal{F} = \{\sigma(x^{*}, \Gamma_{n}(\cdot)): x^{*} \in U^{o}, n \in \mathbb{N}\}$ is uniformly integrable. 
By hypothesis, for each $x^{*} \in U^{o}$ there exists 
$Z^{x^{*}} \in \Sigma$ with $\mu(Z^{x^{*}}) =0$ such that 
for each $n \in \mathbb{N}$, we have
\begin{equation*}
\begin{split}
|\sigma(x^{*}, \Gamma_{n}(t))| \leq g(t)
~\text{ for every }t \in  \Omega \setminus Z^{x^{*}}
\end{split}
\end{equation*}
and, therefore, for any $A \in \Sigma$, 
\begin{equation*}
\begin{split}
|\sigma(x^{*}, \Gamma_{n}(t))| \leq g(t)
~\text{ for every }t \in A \setminus Z^{x^{*}}.
\end{split}
\end{equation*}
Then, 
\begin{equation}\label{eqC_Dominated_Convergence_Theorem.1}
\begin{split}
\int_{A} |\sigma(x^{*}, \Gamma_{n}(t))| d \mu  \leq \int_{A} g(t) d \mu
\qquad
(A \in \Sigma, n \in \mathbb{N}, x^{*} \in U^{o})
\end{split}
\end{equation}
and, in particular, 
\begin{equation*}
\begin{split}
||\sigma(x^{*}, \Gamma_{n})||_{1} = \int_{\Omega} |\sigma(x^{*}, \Gamma_{n}(t))| d \mu  
\leq ||g||_{1}
~\text{ for all }n \in \mathbb{N} \text{ and }x^{*} \in U^{o}.
\end{split}
\end{equation*}

Since $g$ is Lebesgue integrable, given $\varepsilon >0$ there exists $\delta >0$ such that 
for each $A \in \Sigma$, we have
\begin{equation*}
\begin{split}
\mu(A) < \delta \Rightarrow \int_{A} g(t) d \mu < \varepsilon.
\end{split}
\end{equation*}
The last result together with \eqref{eqC_Dominated_Convergence_Theorem.1} yields 
\begin{equation*}
\begin{split}
\mu(A) < \delta 
\Rightarrow 
\sup \left \{\int_{A} |\sigma(x^{*}, \Gamma_{n}(t))| d \mu : x^{*} \in U^{o}, n \in \mathbb{N} \right \} \leq  \int_{A} g(t) d \mu < \varepsilon.
\end{split}
\end{equation*}
This means that the family $\mathcal{F}$ is uniformly integrable, and this completes the proof.
\end{proof}


\section{An application}

There is a close relation between fuzzy maps Pettis integrable and set valued maps Pettis integrable,  
for example the relation between 
Radon-Nikodym derivative of fuzzy number multimeasures 
and 
Radon-Nikodym derivative of multimeasures, 
see \cite{BAN, STOJ_RNT, WU_RNT, ZHOU, KAL_GFNM, KAL_FNM, KAL_FNM_Frechet}.
Moreover, in this section, 
we are going to show that Theorem \ref{T_Pettis_Chr_CWK} holds also 
for fuzzy maps scalarly integrable, see Theorem \ref{T_FPettis_Chr_CWK}. 
We first need to introduce some definitions 
and an auxiliary lemma, Lemma \ref{L_CHR_FP}.

\begin{definition}\label{D_Fuzzy_Number}
A map $u: X \to [0,1]$ is said to be a \textit{fuzzy set} on $X$
(see Zadeh \cite{ZAD}).
For each $\alpha \in (0,1]$, denote by 
\begin{equation*}
\begin{split}
[u]^{\alpha} = \{x \in X: u(x) \geq \alpha\}
\end{split}
\end{equation*}
the $\alpha$-level set of $u$. 
We say that the fuzzy set $u$ on $X$ is a \textit{generalized fuzzy number} on $X$, 
if $[u]^{\alpha} \in bcc(X)$ for all $\alpha \in (0,1]$.  
We denote by $\mathcal{F}_{bcc}(X)$ the set of all generalized fuzzy numbers 
on $X$. 
In particular, 
by $\mathcal{F}_{cwk}(X)$ is denoted the set of all generalized fuzzy numbers 
$u \in \mathcal{F}_{bcc}(X)$ such that 
$[u]^{\alpha} \in cwk(X)$  for all $\alpha \in (0,1]$. 
We write $\widetilde{\theta} = \mathbbm{1}_{\{\theta\}}$, 
where  $\mathbbm{1}_{\{\theta\}}$ is the characteristic function of the set $\{\theta\}$.
\end{definition}
It is easy to see that \cite[Lemma 3.1]{KAL_GFNM} holds also for 
the complete Hausdorff locally convex topological vector space $X$.

\begin{definition}\label{D_FNMultimeasure}
Given a fuzzy map $M : \Sigma \to \mathcal{F}_{cwk}(X)$ we define 
the $\alpha$-cut $M^{\alpha}$ of $M$ as follows 
\begin{equation*}
\begin{split}
M^{\alpha} : \Sigma \to cwk(X),
\quad
M^{\alpha}(A)=[M(A)]^{\alpha} = \{x \in X: M(A)(x) \geq \alpha\}
\quad
(\alpha \in (0,1]).
\end{split}
\end{equation*}
We say that  $M$ is a \textit{fuzzy number multimeasure}, 
if for each  $\alpha \in (0,1]$ $M^{\alpha}$ is a multimeasure.    
\end{definition}

\begin{definition}\label{D_Fuzzy_Integrals}
Given a fuzzy map  $\Gamma : \Omega \to \mathcal{F}_{bcc}(X)$ 
we can define the $\alpha$-cut $\Gamma^{\alpha}$ of $\Gamma$ as follows 
\begin{equation*}
\begin{split}
\Gamma^{\alpha} : \Omega \to bcc(X), 
\quad
\Gamma^{\alpha}(t) = [\Gamma(t)]^{\alpha} = \{x \in X : \Gamma(t)(x) \geq \alpha\}
\quad
(\alpha \in (0,1]).  
\end{split}
\end{equation*}
We say that the fuzzy map 
$\Gamma$ 
is scalarly measurable (scalarly integrable) if 
the map $\sigma(x^{*}, \Gamma^{\alpha}(\cdot))$ is measurable (Lebesgue integrable) 
for each $x^{*} \in X^{*}$ 
and $\alpha \in (0,1]$. 
The fuzzy map $\Gamma$ is said to be \textit{Pettis integrable} in $\mathcal{F}_{bcc}(X)$ ($\mathcal{F}_{cwk}(X)$) if 
\begin{itemize}
\item[(i)]
$\Gamma$ is scalarly integrable, 
\item[(ii)]
for each $A \in \Sigma$ there exists 
a generalized fuzzy number $M(A) \in \mathcal{F}_{bcc}(X)$ ($M(A) \in \mathcal{F}_{cwk}(X)$) such that  
\begin{equation*}
\begin{split}
\sigma \left ( x^{*}, [M(A)]^{\alpha} \right ) 
= 
\int_{A} \sigma(x^{*}, \Gamma^{\alpha}(t)) d \mu
~\text{ for all }x^{*} \in X^{*}, \alpha \in (0,1].
\end{split}
\end{equation*}
\end{itemize}
We call $M(A)$ the \textit{fuzzy Pettis integral} of $\Gamma$ over $A$ and 
set $(FP)\int_{A} \Gamma(t) d \mu = M(A)$. 
It is clear that the map $M$ defined by  $(FP)\int_{A} \Gamma(t) d \mu = M(A)$ for all $A \in \Sigma$, 
is a fuzzy number multimeasure. 
\end{definition}

\begin{lemma}\label{L_CHR_FP}
Let $\Gamma : \Omega \to \mathcal{F}_{cwk}(X)$ be a fuzzy map. 
Then the following statements are equivalent:
\begin{itemize}
\item[(i)]
$\Gamma$ is Pettis integrable in $\mathcal{F}_{cwk}(X)$,
\item[(ii)]
for each $\alpha \in (0,1]$ the $\alpha$-cut $\Gamma^{\alpha}$ of $\Gamma$  
is Pettis integrable.
\end{itemize}
\end{lemma}
\begin{proof}
By Definition \ref{D_Fuzzy_Integrals} follows immediately that $(i) \Rightarrow (ii)$. 
Assume now that $(ii)$ holds and let 
\begin{equation*}
\begin{split}
M_{\alpha}(A): \Sigma \to cwk(X),
\quad
M_{\alpha}(A) 
= (P)\int_{A} \Gamma^{\alpha}(t) d \mu, 
\quad
(\alpha \in (0,1]).
\end{split}
\end{equation*}
It is enough to prove that for each $A \in \Sigma$ 
there exists a fuzzy number $M(A) \in \mathcal{F}_{cwk}(X)$ such that
\begin{equation}\label{eqT_CHR_Br.0}
\begin{split}
[M(A)]^{\alpha} = M_{\alpha}(A)
~\text{ for every }\alpha \in (0,1].
\end{split}
\end{equation}

$\bullet$ 
Given two real numbers $\alpha$ and $\beta$ such that 
$0<\alpha \leq \beta \leq 1$, $\Gamma^{\beta}(t) \subset \Gamma^{\alpha}(t)$ for all $t \in \Omega$. 
Hence,
\begin{equation*} 
\begin{split}
\sigma(x^{*}, M_{\beta}(A))
=
\int_{A} \sigma(x^{*}, \Gamma^{\beta}(t)) d \mu 
\leq 
\int_{A} \sigma(x^{*}, \Gamma^{\alpha}(t)) d \mu 
=
\sigma(x^{*}, M_{\alpha}(A))
\quad
(x^{*} \in X^{*}, A \in \Sigma).
\end{split}
\end{equation*}
Thanks to  the Hahn-Banach separation
theorem (\cite[Theorem V.2.10]{DUN}) the last result yields that $M_{\beta}(A) \subset M_{\alpha}(A)$ 
for all $A \in \Sigma$.

$\bullet$
Fix an arbitrary $A \in \Sigma$ and $\alpha \in (0,1]$,  
and let $(\alpha_{n})$ be an increasing sequence with terms in $(0,1]$ 
converging to $\alpha$. 
Then,
\begin{equation}\label{eqT_CHR_Br.1}
\begin{split}
M_{\alpha}(A) \subset \bigcap_{n} M_{\alpha_{n}}(A)
\end{split}
\end{equation}
To prove the reverse inclusion, 
let $x_{0} \in \bigcap_{n} M_{\alpha_{n}}(A) \setminus M_{\alpha}(A)$. 
Then, 
by the Hahn-Banach separation
theorem (\cite[Theorem V.2.10]{DUN}) there exists $x^{*}_{0} \in X^{*}$, $r \in \mathbb{R}$ 
and $\delta >0$ such that 
\begin{equation*}
\begin{split}
\sigma(x^{*}_{0}, M_{\alpha}(A)) \leq r < r + \delta \leq x^{*}_{0}(x_{0})
\end{split}
\end{equation*}
and since  
$x_{0} \in M_{\alpha_{n}}(A)$ for every $n \in \mathbb{N}$, 
it follows that 
\begin{equation}\label{eqT_CHR_Br.2}
\begin{split}
\sigma(x^{*}_{0}, M_{\alpha}(A)) \leq r < r + \delta 
\leq \lim_{n \to \infty} \sigma(x^{*}_{0}, M_{\alpha_{n}}(A)).
\end{split}
\end{equation}
On the other hand, 
from the fact  
$\Gamma^{\alpha}(t) = \bigcap_{n} \Gamma^{\alpha_{n}}(t)$ for every $t \in \Omega$ 
(see \cite[Lemma 3.1]{KAL_GFNM}), 
it follows that 
\begin{equation*}
\begin{split}
\lim_{n \to \infty} \sigma(x^{*}, \Gamma^{\alpha_{n}}(t)) = \sigma(x^{*}, \Gamma^{\alpha}(t))
\quad
(x^{*} \in X^{*}, t \in \Omega), 
\end{split}
\end{equation*}
and since 
\begin{equation*}
\begin{split}
|\sigma(x^{*}, \Gamma^{\alpha_{n}}(t))|
\leq
\max 
\{
|\sigma(x^{*}, \Gamma^{\alpha_{1}}(t))|, |\sigma(x^{*}, \Gamma^{\alpha}(t))|
\}
\quad
(x^{*} \in X^{*}, t \in \Omega, n \in \mathbb{N})
\end{split}
\end{equation*}
we obtain by Dominated Convergence Theorem that 
\begin{equation*}
\begin{split}
\lim_{n \to \infty} \sigma(x^{*}, M_{\alpha_{n}}(A))
=&
\lim_{n \to \infty} \int_{A} \sigma(x^{*}, \Gamma^{\alpha_{n}}(t)) d \mu 
\\
=& 
\int_{A} \sigma(x^{*}, \Gamma^{\alpha}(t)) d \mu = \sigma(x^{*}, M_{\alpha}(A))
\quad
(x^{*} \in X^{*}).  
\end{split}
\end{equation*}
The last result together with \eqref{eqT_CHR_Br.2} yields the contradict 
$\sigma(x^{*}_{0}, M_{\alpha}(A)) < \sigma(x^{*}_{0}, M_{\alpha}(A))$.
Hence, 
$x \not\in M_{\alpha}(A)$ implies that $x \not\in \bigcap_{n} M_{\alpha_{n}}(A)$. 
Thus, $\bigcap_{n} M_{\alpha_{n}}(A) \subset M_{\alpha}(A)$ 
and, consequently by \eqref{eqT_CHR_Br.1} follows that 
$\bigcap_{n} M_{\alpha_{n}}(A) = M_{\alpha}(A)$.

$\bullet$ 
Thus, we  have that for each $A \in \Sigma$ the family 
$\{M_{\alpha}(A) : \alpha \in (0,1]\}$ satisfies $(i), (ii)$ and $(iii)$ in \cite[Lemma 3.1]{KAL_GFNM}.
Therefore, for each $A \in \Sigma$, 
there is a unique fuzzy number $M(A) \in \mathcal{F}_{cwk}(X)$ such that 
\eqref{eqT_CHR_Br.0} holds. 
Consequently, the fuzzy map $\Gamma$ is  Pettis integrable with 
$M(A) = (FP)\int_{A} \Gamma(t) d \mu$ for all $A \in \Sigma$, 
and this completes the proof.
\end{proof}


Given a fuzzy map $\Gamma: \Omega \to \mathcal{F}_{cwk}(X)$ scalarly integrable define
\begin{equation}\label{eq_FVAR_Gamma}
\begin{split}
\varphi_{\Gamma} : X^{*} \to L = \prod_{\alpha \in (0,1]} L_{1}(\mu), 
\quad
\varphi_{\Gamma}(x^{*}) = (\sigma(x^{*}, \Gamma^{\alpha}(t))_{\alpha \in (0,1]}, 
\end{split}
\end{equation}
where $L$ is equipped with the product topology $\mathfrak{T}_{P}$ 
and $\sigma(L, L^{*})$ is the weak topology in $L$. 


We are now ready to prove the main result of this section.

\begin{theorem}\label{T_FPettis_Chr_CWK}
Let $\Gamma : \Omega \to \mathcal{F}_{cwk}(X)$ be a scalarly integrable fuzzy map. 
Then the following statements are equivalent:
\begin{itemize}
\item[(i)]
$\Gamma$ is Pettis integrable in $\mathcal{F}_{cwk}(X)$, 
\item[(ii)]
$\varphi_{\Gamma}$ is $\tau(X^{*},X)$-$\mathfrak{T}_{P}$ continuous,
\item[(iii)]
$\varphi_{\Gamma}$ is $\tau(X^{*},X)$-$\sigma(L, L^{*})$ continuous,
\item[(iv)]
for each $\alpha \in (0,1]$ and for each $A \in \Sigma$ the map
\begin{equation*}
\begin{split}
\varphi^{A}_{\alpha} : X^{*} \to \mathbb{R},
\quad
\varphi^{A}_{\alpha}(x^{*}) = \int_{A} \sigma(x^{*}, \Gamma^{\alpha}(t)) d \mu
\end{split}
\end{equation*}
is  $\tau(X^{*},X)$-continuous. 
\end{itemize}
\end{theorem}
\begin{proof}
By \eqref{eq_FVAR_Gamma}, we obtain
\begin{equation*}
\begin{split}
\varphi_{\Gamma^{\alpha}} = \pi_{\alpha} \circ \varphi_{\Gamma}
~\text{ for every } \alpha \in (0,1], 
\end{split}
\end{equation*}
where $\pi_{\alpha}$ is the $\alpha$-projection and 
$\varphi_{\Gamma^{\alpha}}$ is the sublinear map defined by $\Gamma^{\alpha}$, 
see Definition \ref{D_V_Gamma}.
According to \cite[Theorem 4.3(1), p.137]{SCH}, we have
\begin{equation*}
\begin{split}
\sigma(L, L^{*}) = \prod_{\alpha \in (0,1]} \sigma( L_{1}, L_{1}^{*}), 
\end{split}
\end{equation*}
where $\sigma( L_{1}, L_{1}^{*})$ is the weak topology in $L_{1}(\mu)$.

Thanks to \cite[5.2, p.51]{SCH}, we have
\begin{itemize}
\item[(a)]
$\varphi_{\Gamma}$ is $\tau(X^{*},X)$-$\mathfrak{T}_{P}$-continuous, 
if and only if for each $\alpha \in (0,1]$, the map
$\varphi_{\Gamma^{\alpha}}$ is $\tau(X^{*},X)$-$||\cdot||_{1}$-continuous,
\item[(b)]
$\varphi_{\Gamma}$ is $\tau(X^{*},X)$-$\sigma(L, L^{*})$-continuous, 
if and only if for each $\alpha \in (0,1]$, the map $\varphi_{\Gamma^{\alpha}}$ is 
$\tau(X^{*},X)$-$\sigma( L_{1}, L_{1}^{*})$-continuous. 
\end{itemize}
Hence, Lemma \ref{L_CHR_FP} together with Theorem \ref{T_Pettis_Chr_CWK} yields that 
$(i) \Leftrightarrow (ii) \Leftrightarrow (iii)$.

It remains to prove that $(iii) \Rightarrow (iv) \Rightarrow (i)$. 
Assume that $(iii)$ holds. 
Then, for each $\alpha \in (0,1]$, the map $\varphi_{\Gamma^{\alpha}}$ is 
$\tau(X^{*},X)$-$\sigma( L_{1}, L_{1}^{*})$ continuous (see $(b)$). 
Hence, by Theorem \ref{T_Pettis_Chr_CWK} follows that $(iv)$ holds.

We now assume that $(iv)$ holds. 
Then, we obtain by Theorem \ref{T_Pettis_Chr_CWK} that for each $\alpha \in (0,1]$, 
the $\alpha$-cut $\Gamma^{\alpha}$ of $\Gamma$ is Pettis integrable. 
Therefore, by Lemma \ref{L_CHR_FP} the fuzzy map $\Gamma$ is Pettis integrable, and this completes the proof.
\end{proof}

\end{document}